\documentclass{amsart}
\usepackage{amsthm}
\usepackage{amssymb}

\makeatletter

\theoremstyle{plain}
\newtheorem{thm}{\protect\theoremname}
  \theoremstyle{plain}
  \newtheorem{cor}[thm]{\protect\corollaryname}
  \theoremstyle{remark}
  \newtheorem*{rem*}{\protect\remarkname}

\makeatother

  \providecommand{\corollaryname}{Corollary}
  \providecommand{\remarkname}{Remark}
\providecommand{\theoremname}{Theorem}

\begin{document}

\title{Identities of symmetry for higher-order $q$-Bernoulli polynomials}

\author{Dae San Kim \and Taekyun Kim}
\begin{abstract}
Recently, the higher-order Carlitz's $q$-Bernoulli polynomials are
represented as $q$-Volkenborn integral on $\mathbb{Z}_{p}$  by Kim. 
A question was asked in \cite{key-13} as to finding the extended formulae of
symmetries for Bernoulli polynomials which are related to Carlitz
$q$-Bernoulli polynomials. In this paper, we give some
new identities of symmetry for the higher-order Carlitz's $q$-Bernoulli
polynomials which are derived from multivariate $q$-Volkenborn integrals
on $\mathbb{Z}_{p}$. We note that they are a partial answer to that question.
\end{abstract}



\maketitle

\section{Introduction}

$\,$\global\long\def\Zp{\mathbb{Z}_{p}}

Let $p$ be a fixed prime number. Throughout this paper, $\Zp$, $\mathbb{Q}_{p}$
and $\mathbb{C}_{p}$ will, respectively, denote the ring of $p$-adic
rational integers, the field of $p$-adic rational numbers and the
completion of algebraic closure of $\mathbb{Q}_{p}$. The $p$-adic
absolute value in $\mathbb{C}_{p}$ is normalized so that $\left|p\right|_{p}=p^{-1}$.
Let $q$ be an indeterminate in $\mathbb{C}_{p}$ with $\left|1-q\right|_{p}<p^{-\frac{1}{p-1}}$.
We say that $f$ is uniformly differentiable function at a point $a\in\Zp$,
if the difference quotient, 
\[
F_{f}\,:\,\Zp\times\Zp\rightarrow\Zp\quad\mbox{by }F_{f}\left(x,y\right)=\frac{f\left(x\right)-f\left(y\right)}{x-y},
\]
have a limit $l=f^{\prime}\left(a\right)$ as $\left(x,y\right)\rightarrow\left(a,a\right)$. If 
$f$ is uniformly differentiable on $\Zp$, we denote this property by $f\in UD\left(\Zp\right)$.

For $f\in UD\left(\Zp\right)$ , the $q$-Volkenborn integral is defined
by Kim to be 
\begin{equation}
I_{q}\left(f\right)=\int_{\Zp}f\left(x\right)d\mu_{q}\left(x\right)=\lim_{N\rightarrow\infty}\frac{1}{\left[p^{N}\right]_{q}}\sum_{x=0}^{p^{N}-1}f\left(x\right)q^{x},\label{eq:1}
\end{equation}
where $\left[x\right]_{q}=\frac{1-q^{x}}{1-q}$, (see \cite{key-14,key-13,key-12}).

From (\ref{eq:1}), we note that 
\begin{equation}
qI_{q}\left(f_{1}\right)=I_{q}\left(f\right)+\left(q-1\right)f\left(0\right)+\frac{q-1}{\log q}f^{\prime}\left(0\right)\label{eq:2}
\end{equation}
where $f_{1}\left(x\right)=f\left(x+1\right)$.

As is well known, the Bernoulli polynomials are defined by the generating
function to be
\begin{equation}
\frac{t}{e^{t}-1}e^{xt}=e^{B\left(x\right)t}=\sum_{n=0}^{\infty}B_{n}\left(x\right)\frac{t^{n}}{n!}.\label{eq:3}
\end{equation}

When $x=0$, $B_{n}=B_{n}\left(0\right)$ are called the Bernoulli
numbers.

By (\ref{eq:3}), we get 
\begin{equation}
\left(B+1\right)^{n}-B_{n}=\begin{cases}
1 & \mbox{if }n=1\\
0 & \mbox{if }n>1,
\end{cases}\:\mbox{and }B_{0}=1.\label{eq:4}
\end{equation}

In \cite{key-3}, Carlitz defined $q$-Bernoulli numbers as follows
: 
\begin{equation}
\beta_{0,q}=1,\qquad q\left(q\beta_{q}+1\right)^{n}-\beta_{n,q}=\begin{cases}
1, & \mbox{if }n=1\\
0, & \mbox{if }n>1,
\end{cases}\label{eq:5}
\end{equation}
with the usual convention about replacing $\beta_{q}^{i}$ by $\beta_{i,q}$.

From (\ref{eq:4}) and (\ref{eq:5}), we note that ${\displaystyle \lim_{q\rightarrow1}\beta_{n,q}=B_{n}.}$ 

The $q$-Bernoulli polynomials are given by 
\begin{align}
\beta_{n,q}\left(x\right) & =\sum_{l=0}^{n}\dbinom{n}{l}q^{lx}\beta_{l,q}\left[x\right]_{q}^{n-l}\label{eq:6}\\
 & =\left(q^{x}\beta_{q}+\left[x\right]_{q}\right)^{n},\quad\left(n\ge0\right),\quad\mbox{(see \cite{key-2,key-3,key-4,key-12}}).\nonumber 
\end{align}

In \cite{key-14}, Kim proved that Carlitz $q$-Bernoulli polynomials
can be written by $q$-Volkenborn integral on $\Zp$ as follows :
\begin{align}
\beta_{n,q}\left(x\right) & =\int_{\Zp}\left[x+y\right]_{q}^{n}d\mu_{q}\left(y\right)\label{eq:7}\\
 & =\sum_{l=0}^{n}\dbinom{n}{l}\left[x\right]_{q}^{n-l}q^{lx}\int_{\Zp}\left[y\right]_{q}^{l}d\mu_{q}\left(x\right).\nonumber 
\end{align}

Thus, by (\ref{eq:7}), we get 
\begin{equation}
\beta_{n,q}=\int_{\Zp}\left[x\right]_{q}^{n}d\mu_{q}\left(x\right),\quad\left(n\ge0\right).\label{eq:8}
\end{equation}

From (\ref{eq:2}), we note that 
\begin{equation}
q\int_{\Zp}\left[x+1\right]_{q}^{n}d\mu_{q}\left(x\right)-\int_{\Zp}\left[x\right]_{q}^{n}d\mu_{q}\left(x\right)=\begin{cases}
q-1 & \mbox{if }n=0\\
1 & \mbox{if }n=1\\
0 & \mbox{if }n>1.
\end{cases}\label{eq:9}
\end{equation}

By (\ref{eq:7}), (\ref{eq:8}) and (\ref{eq:9}), we see that 
\begin{equation}
\beta_{0,q}=1,\quad q\left(q\beta_{q}+1\right)^{n}-\beta_{n,q}=\begin{cases}
1 & \mbox{if }n=1\\
0 & \mbox{if }n>1.
\end{cases}\label{eq:10}
\end{equation}

Let 
\begin{align}
I_{1}\left(f\right)=\lim_{q\rightarrow1}I_{q}\left(f\right) & =\int_{\Zp}f\left(x\right)d\mu_{1}\left(x\right)\label{eq:11}\\
 & =\lim_{N\rightarrow\infty}\frac{1}{p^{N}}\sum_{x=0}^{p^{N}-1}f\left(x\right),\quad\left(\mbox{see \cite{key-14,key-23,key-24}}\right).\nonumber 
\end{align}

Then, by (\ref{eq:2}), we get 
\begin{equation}
I_{1}\left(f_{1}\right)-I_{1}\left(f\right)=f^{\prime}\left(0\right).\label{eq:12}
\end{equation}

Let us take $f\left(x\right)=e^{tx}$. Then we have 
\begin{equation}
\int_{\Zp}e^{xt}d\mu_{1}\left(x\right)=\frac{t}{e^{t}-1}=\sum_{n=0}^{\infty}B_{n}\frac{t^{n}}{n!},\label{eq:13}
\end{equation}
and 
\begin{equation}
\int_{\Zp}e^{\left(x+y\right)t}d\mu_{1}\left(y\right)=\left(\frac{t}{e^{t}-1}\right)e^{xt}=\sum_{n=0}^{\infty}B_{n}\left(x\right)\frac{t^{n}}{n!}.\label{eq:14}
\end{equation}

For $r\in\mathbb{N}$, the higher-order Bernoulli polynomials are
defined by the generating function to be 
\begin{equation}
\left(\frac{t}{e^{t}-1}\right)^{r}e^{xt}=\underset{r-\mbox{times}}{\underbrace{\left(\frac{t}{e^{t}-1}\right)\times\cdots\times\left(\frac{t}{e^{t}-1}\right)}}e^{xt}=\sum_{n=0}^{\infty}B_{n}^{\left(r\right)}\left(x\right)\frac{t^{n}}{n!}.\label{eq:15}
\end{equation}

By (\ref{eq:14}), we get 
\begin{align}
\int_{\Zp}\cdots\int_{\Zp}e^{\left(x+y_{1}+\cdots+y_{r}\right)}d\mu_{1}\left(y_{1}\right)\cdots d\mu_{1}\left(y_{r}\right) & =\left(\frac{t}{e^{t}-1}\right)^{r}e^{xt}\label{eq:16}\\
 & =\sum_{n=0}^{\infty}B_{n}^{\left(r\right)}\left(x\right)\frac{t^{n}}{n!}.\nonumber 
\end{align}

In \cite{key-3,key-4}, Carlitz introduced the $q$-extension of higher-order
Bernoulli polynomials as follows : 
\begin{equation}
\beta_{n,q}^{\left(r\right)}\left(x\right)=\frac{1}{\left(1-q\right)^{n}}\sum_{l=0}^{n}\dbinom{n}{l}\left(-1\right)^{l}q^{lx}\left(\frac{l+1}{\left[l+1\right]_{q}}\right)^{r},\label{eq:17}
\end{equation}
where $n\ge0$ and $r\in\mathbb{N}$.

Note that ${\displaystyle \lim_{q\rightarrow1}}\beta_{n,q}^{\left(r\right)}\left(x\right)=B_{n}^{\left(r\right)}\left(x\right)$. 

From (\ref{eq:16}), we note that 
\begin{equation}
\int_{\Zp}\cdots\int_{\Zp}\left(x+y_{1}+\cdots+y_{r}\right)^{n}d\mu_{1}\left(y_{1}\right)\cdots d\mu_{1}\left(y_{r}\right)=B_{n}^{\left(r\right)}\left(x\right),\label{eq:18}
\end{equation}
where $n\ge0$ and $r\in\mathbb{N}$.

In this paper, we consider $q$-extensions of (\ref{eq:17}) which
are related to higher-order Carlitz's $q$-Bernoulli polynomials.
The purpose of this paper is to give some new and interesting identities
of symmetry for the higher-order Carlitz's $q$-Bernoulli polynomials
which are derived from multivariate $q$-Volkenborn integral on $\Zp$.

\section{Identities of symmetry for higher-order $q$-Bernoulli polynomials}

$\,$

In the sense of $q$-extension of (\ref{eq:18}), we observe the following
equation (\ref{eq:19}) 
\begin{align}
 & \int_{\Zp}\cdots\int_{\Zp}\left[x+y_{1}+\cdots+y_{r}\right]_{q}^{n}d\mu_{q}\left(y_{1}\right)\cdots d\mu_{q}\left(y_{n}\right)\label{eq:19}\\
= & \frac{1}{\left(1-q\right)^{n}}\sum_{l=0}^{n}\dbinom{n}{l}\left(-1\right)^{l}q^{lx}\left(\frac{l+1}{\left[l+1\right]_{q}}\right)^{r}.\nonumber 
\end{align}

Thus, by (\ref{eq:17}) and (\ref{eq:19}), we get 
\begin{equation}
\beta_{n,q}^{\left(r\right)}\left(x\right)=\int_{\Zp}\cdots\int_{\Zp}\left[x+y_{1}+\cdots+y_{r}\right]_{q}^{n}d\mu_{q}\left(y_{1}\right)\cdots d\mu_{q}\left(y_{r}\right),\label{eq:20}
\end{equation}
where $n\ge0$ and $r\in\mathbb{N}$.

Let us consider the generating function of $\beta_{n,q}^{\left(r\right)}\left(x\right)$
as follows : 
\begin{equation}
\sum_{n=0}^{\infty}\beta_{n,q}^{\left(r\right)}\left(x\right)\frac{t^{n}}{n!}=\int_{\Zp}\cdots\int_{\Zp}e^{\left[x+y_{1}+\cdots+y_{r}\right]_{q}t}d\mu_{q}\left(y_{1}\right)\cdots d\mu_{q}\left(y_{r}\right).\label{eq:21}
\end{equation}

For $w_{1},\, w_{2}\in\mathbb{N}$, we have
\begin{align}
 & \frac{1}{\left[w_{1}\right]_{q}^{r}}\int_{\Zp}\cdots\int_{\Zp}e^{\left[w_{1}w_{2}x+w_{2}\sum_{l=1}^{r}j_{l}+w_{1}\sum_{l=1}^{r}y_{l}\right]_{q}t}d\mu_{q^{w_{1}}}\left(y_{1}\right)\cdots d\mu_{q^{w_{1}}}\left(y_{r}\right)\label{eq:22}\\
= & \lim_{N\rightarrow\infty}\left(\frac{1}{\left[w_{1}\right]_{q}\left[p^{N}\right]_{q^{w_{1}}}}\right)^{r}\sum_{y_{1},\cdots,\, y_{r}=0}^{p^{N}-1}e^{\left[w_{1}w_{2}x+w_{2}\sum_{l=1}^{r}j_{l}+w_{1}\sum_{l=1}^{r}y_{l}\right]_{q}t}q^{w_{1}\left(y_{1}+\cdots+y_{r}\right)}\nonumber \\
= & \lim_{N\rightarrow\infty}\left(\frac{1}{\left[w_{1}\right]_{q}\left[w_{2}p^{N}\right]_{q^{w_{1}}}}\right)^{r}\nonumber \\
 & \times\sum_{y_{1},\cdots,y_{r}=0}^{w_{2}p^{N}-1}e^{\left[w_{1}w_{2}x+w_{2}\sum_{l=1}^{r}j_{l}+w_{1}\sum_{l=1}^{r}y_{l}\right]_{q}t}q^{w_{1}\left(y_{1}+\cdots+y_{r}\right)}\nonumber \\
= & \lim_{N\rightarrow\infty}\left(\frac{1}{\left[w_{1}w_{2}p^{N}\right]_{q}}\right)^{r}\nonumber \\
 & \times\sum_{i_{1},\cdots,\, i_{r}=0}^{w_{2}-1}\sum_{y_{1},\cdots,y_{r}=0}^{p^{N}-1}e^{\left[w_{1}w_{2}x+\sum_{l=1}^{r}\left(w_{2}j_{l}+w_{1}i_{l}+w_{1}w_{2}y_{l}\right)\right]_{q}t}q^{w_{1}\sum_{l=1}^{r}\left(i_{l}+w_{2}y_{l}\right)}.\nonumber 
\end{align}

Thus, by (\ref{eq:22}), we get 
\begin{align}
 & \frac{1}{\left[w_{1}\right]_{q}^{r}}\sum_{j_{1},\cdots,j_{r}=0}^{w_{1}-1}q^{w_{2}\sum_{l=1}^{r}j_{l}}\label{eq:23}\\
 & \times\int_{\Zp}\cdots\int_{\Zp}e^{\left[w_{1}w_{2}x+\sum_{l=1}^{r}\left(w_{2}j_{l}+w_{1}y_{l}\right)\right]_{q}t}d\mu_{q^{w_{1}}}\left(y_{1}\right)\cdots d\mu_{q^{w_{1}}}\left(y_{r}\right)\nonumber \\
= & \lim_{N\rightarrow\infty}\left(\frac{1}{\left[w_{1}w_{2}p^{N}\right]_{q}}\right)^{r}\sum_{j_{1},\cdots,\, j_{r}=0}^{w_{1}-1}\sum_{i_{1},\cdots,\, i_{r}=0}^{w_{2}-1}\nonumber \\
 & \times\sum_{y_{1},\cdots,\, y_{r}=0}^{p^{N}-1}e^{\left[w_{1}w_{2}x+\sum_{l=1}^{r}\left(w_{2}j_{l}+w_{1}i_{l}+w_{1}w_{2}y_{l}\right)\right]_{q}t}q^{\sum_{l=1}^{r}\left(w_{1}i_{l}+w_{2}j_{l}+w_{1}w_{2}y_{l}\right)}.\nonumber 
\end{align}

By the same method as (\ref{eq:23}), we get 
\begin{align}
 & \frac{1}{\left[w_{2}\right]_{q}^{r}}\sum_{j_{1},\cdots,j_{r}=0}^{w_{2}-1}q^{w_{1}\sum_{l=1}^{r}j_{l}}\label{eq:24}\\
 & \times\int_{\Zp}\cdots\int_{\Zp}e^{\left[w_{1}w_{2}x+\sum_{l=1}^{r}\left(w_{1}j_{l}+w_{2}y_{l}\right)\right]_{q}t}d\mu_{q^{w_{2}}}\left(y_{1}\right)\cdots d\mu_{q^{w_{2}}}\left(y_{l}\right)\nonumber \\
= & \lim_{N\rightarrow\infty}\left(\frac{1}{\left[w_{1}w_{2}p^{N}\right]_{q}}\right)^{r}\sum_{j_{1},\cdots,j_{r}=0}^{w_{2}-1}\sum_{i_{1},\cdots,i_{r}=0}^{w_{1}-1}\nonumber \\
 & \times\sum_{y_{1},\cdots,y_{r}=0}^{p^{N}-1}e^{\left[w_{1}w_{2}x+\sum_{l=1}^{r}\left(w_{1}j_{l}+w_{2}i_{l}+w_{1}w_{2}y_{l}\right)\right]_{q}t}q^{\sum_{l=1}^{r}\left(w_{2}i_{l}+w_{1}j_{l}+w_{1}w_{2}y_{l}\right)}.\nonumber 
\end{align}

Therefore, by (\ref{eq:23}), we obtain the following theorem.
\begin{thm}
\label{thm:1} For $w_{1},\, w_{2}\in\mathbb{N}$, we have
\begin{align*}
 & \frac{1}{\left[w_{1}\right]_{q}^{r}}\sum_{j_{1},\cdots,j_{r}=0}^{w_{1}-1}q^{w_{2}\sum_{l=1}^{r}j_{l}}\\
 & \times\int_{\Zp}\cdots\int_{\Zp}e^{\left[w_{1}w_{2}x+\sum_{l=1}^{r}\left(w_{2}j_{l}+w_{1}y_{l}\right)\right]_{q}t}d\mu_{q^{w_{1}}}\left(y_{1}\right)\cdots d\mu_{q^{w_{1}}}\left(y_{r}\right)\\
= & \frac{1}{\left[w_{2}\right]_{q}^{r}}\sum_{j_{1},\cdots,j_{r}=0}^{w_{2}-1}q^{w_{1}\sum_{l=1}^{r}j_{l}}\\
 & \times\int_{\Zp}\cdots\int_{\Zp}e^{\left[w_{1}w_{2}x+\sum_{l=1}^{r}\left(w_{1}j_{l}+w_{2}y_{l}\right)\right]_{q}t}d\mu_{q^{w_{2}}}\left(y_{1}\right)\cdots d\mu_{q^{w_{2}}}\left(y_{r}\right).
\end{align*}

\end{thm}
$\,$

It is easy to show that 
\begin{align}
 & \left[w_{1}w_{2}x+w_{2}\left(j_{1}+\cdots+j_{r}\right)+w_{1}\left(y_{1}+\cdots+y_{r}\right)\right]_{q}\label{eq:25}\\
= & \left[w_{1}\right]_{q}\left[w_{2}x+\frac{w_{2}}{w_{1}}\left(j_{1}+\cdots+j_{r}\right)+\left(y_{1}+\cdots+y_{r}\right)\right]_{q^{w_{1}}}.\nonumber 
\end{align}

Therefore, by (\ref{eq:20}), Theorem \ref{thm:1} and (\ref{eq:25}),
we obtain the following corollary, and theorem.
\begin{cor}
\label{cor:2} For $n\ge0$ and $w_{1},\, w_{2}\in\mathbb{N}$, we
have 
\begin{align*}
 & \left[w_{1}\right]_{q}^{n-r}\sum_{j_{1},\cdots,j_{r}=0}^{w_{1}-1}q^{w_{2}\sum_{l=1}^{r}j_{l}}\\
 & \times\int_{\Zp}\cdots\int_{\Zp}\left[w_{2}x+\frac{w_{2}}{w_{1}}\left(j_{1}+\cdots+j_{r}\right)+y_{1}+\cdots+y_{r}\right]_{q^{w_{1}}}^{n}d\mu_{q^{w_{1}}}\left(y_{1}\right)\cdots d\mu_{q^{w_{1}}}\left(y_{r}\right)\\
= & \left[w_{2}\right]_{q}^{n-r}\sum_{j_{1},\cdots,j_{r}=0}^{w_{2}-1}q^{w_{1}\sum_{l=1}^{r}j_{l}}\\
 & \times\int_{\Zp}\cdots\int_{\Zp}\left[w_{1}x+\frac{w_{1}}{w_{2}}\left(j_{1}+\cdots+j_{r}\right)+y_{1}+\cdots+y_{r}\right]_{q^{w_{2}}}^{n}d\mu_{q^{w_{2}}}\left(y_{1}\right)\cdots d\mu_{q^{w_{2}}}\left(y_{r}\right).
\end{align*}

\end{cor}
$\,$
\begin{thm}
\label{thm:3} For $n\ge0$ and $w_{1},\, w_{2}\in\mathbb{N}$, we
have
\begin{align*}
 & \left[w_{1}\right]_{q}^{n-r}\sum_{j_{1},\cdots,j_{r}=0}^{w_{1}-1}q^{w_{2}\left(j_{1}+\cdots+j_{r}\right)}\beta_{n,q^{w_{1}}}^{\left(r\right)}\left(w_{2}x+\frac{w_{2}}{w_{1}}\left(j_{1}+\cdots+j_{r}\right)\right)\\
= & \left[w_{2}\right]_{q}^{n-r}\sum_{j_{1},\cdots,j_{r}=0}^{w_{2}-1}q^{w_{1}\left(j_{1}+\cdots+j_{r}\right)}\beta_{n,q^{w_{2}}}^{\left(r\right)}\left(w_{1}x+\frac{w_{1}}{w_{2}}\left(j_{1}+\cdots+j_{r}\right)\right).
\end{align*}
\end{thm}
\begin{rem*}
Let $w_{2}=1$. Then we have 
\[
\beta_{n,q}^{\left(r\right)}\left(w_{1}x\right)=\left[w_{1}\right]_{q}^{n-r}\sum_{j_{1},\cdots,j_{r}=0}^{w_{1}-1}q^{j_{1}+\cdots+j_{r}}\beta_{n,q^{w_{1}}}^{\left(r\right)}\left(x+\frac{j_{1}+\cdots+j_{r}}{w_{1}}\right).
\]

\end{rem*}
$\,$

By (\ref{eq:20}), we see that 

\begin{align}
 & \int_{\Zp}\cdots\int_{\Zp}\left[w_{2}x+\frac{w_{2}}{w_{1}}\left(j_{1}+\cdots+j_{r}\right)+y_{1}+\cdots+y_{r}\right]_{q^{w_{1}}}^{n}d\mu_{q^{w_{1}}}\left(y_{1}\right)\cdots d\mu_{q^{w_{1}}}\left(y_{r}\right)\label{eq:26}\\
= & \sum_{i=0}^{n}\dbinom{n}{i}\left(\frac{\left[w_{2}\right]_{q}}{\left[w_{1}\right]_{q}}\right)^{i}\left[j_{1}+\cdots+j_{r}\right]_{q^{w_{2}}}^{i}q^{w_{2}\left(n-i\right)\sum_{l=1}^{r}j_{l}}\nonumber \\
 & \times\int_{\Zp}\cdots\int_{\Zp}\left[w_{2}x+\sum_{l=1}^{r}y_{l}\right]_{q^{w_{1}}}^{n-i}d\mu_{q^{w_{1}}}\left(y_{1}\right)\cdots d\mu_{q^{w_{1}}}\left(y_{r}\right)\nonumber \\
= & \sum_{i=0}^{n}\dbinom{n}{i}\left(\frac{\left[w_{2}\right]_{q}}{\left[w_{1}\right]_{q}}\right)^{i}\left[j_{1}+\cdots+j_{r}\right]_{q^{w_{2}}}^{i}q^{w_{2}\left(n-i\right)\sum_{l=1}^{r}j_{l}}\beta_{n-i,q^{w_{1}}}^{\left(r\right)}\left(w_{2}x\right).\nonumber 
\end{align}

From Corollary \ref{cor:2} and (\ref{eq:26}), we have 
\begin{align}
 & \left[w_{1}\right]_{q}^{n-r}\sum_{j_{1},\cdots,j_{r}=0}^{w_{1}-1}q^{w_{2}\sum_{l=1}^{r}j_{l}}\label{eq:27}\\
 & \times\int_{\Zp}\cdots\int_{\Zp}\left[w_{2}x+\frac{w_{2}}{w_{1}}\sum_{l=1}^{r}j_{l}+\sum_{l=1}^{r}y_{l}\right]_{q^{w_{1}}}^{n}d\mu_{q^{w_{1}}}\left(y_{1}\right)\cdots d\mu_{q^{w_{1}}}\left(y_{r}\right)\nonumber \\
= & \sum_{j_{1},\cdots,\, j_{r}=0}^{w_{1}-1}q^{w_{2}\sum_{l=1}^{r}j_{l}}\sum_{i=0}^{n}\dbinom{n}{i}\left[w_{2}\right]_{q}^{i}\left[w_{1}\right]_{q}^{n-i-r}\nonumber \\
 & \times\left[j_{1}+\cdots+j_{r}\right]_{q^{w_{2}}}^{i}q^{w_{2}\left(n-i\right)\sum_{l=1}^{r}j_{l}}\beta_{n-i,q^{w_{1}}}^{\left(r\right)}\left(w_{2}x\right)\nonumber \\
= & \sum_{i=0}^{n}\dbinom{n}{i}\left[w_{2}\right]_{q}^{i}\left[w_{1}\right]_{q}^{n-i-r}\beta_{n-i,q^{w_{1}}}^{\left(r\right)}\left(w_{2}x\right)\nonumber \\
 & \times\sum_{j_{1},\cdots,j_{r}=0}^{w_{1}-1}\left[j_{1}+\cdots+j_{r}\right]_{q^{w_{2}}}^{i}q^{w_{2}\left(n-i+1\right)\sum_{l=1}^{r}j_{l}}\nonumber \\
= & \sum_{i=0}^{n}\dbinom{n}{i}\left[w_{2}\right]_{q}^{n-i}\left[w_{1}\right]_{q}^{i-r}\beta_{i,q^{w_{1}}}^{\left(r\right)}\left(w_{2}x\right)\nonumber \\
 & \times\sum_{j_{1},\cdots,\, j_{r}=0}^{w_{1}-1}\left[j_{1}+\cdots+j_{r}\right]_{q^{w_{2}}}^{n-i}q^{w_{2}\left(i+1\right)\sum_{l=1}^{r}j_{l}}\nonumber \\
= & \sum_{i=0}^{n}\dbinom{n}{i}\left[w_{2}\right]_{q}^{n-i}\left[w_{1}\right]_{q}^{i-r}\beta_{i,q^{w_{1}}}^{\left(r\right)}\left(w_{2}x\right)T_{n,i}^{\left(r\right)}\left(w_{1}|q^{w_{2}}\right),\nonumber 
\end{align}
where 
\begin{equation}
T_{n,i}^{\left(r\right)}\left(w|q\right)=\sum_{j_{1},\cdots,j_{r}=0}^{w-1}\left[j_{1}+\cdots+j_{r}\right]_{q}^{n-i}q^{\left(i+1\right)\sum_{l=1}^{r}j_{l}}.\label{eq:28}
\end{equation}

By the same method as (\ref{eq:28}), we get 
\begin{align}
 & \left[w_{2}\right]_{q}^{n-r}\sum_{j_{1},\cdots,j_{r}=0}^{w_{2}-1}q^{w_{1}\sum_{l=1}^{r}j_{l}}\label{eq:29}\\
 & \times\int_{\Zp}\cdots\int_{\Zp}\left[w_{1}x+\frac{w_{1}}{w_{2}}\sum_{l=1}^{r}j_{l}+\sum_{l=1}^{r}y_{l}\right]_{q^{w_{2}}}^{n}d\mu_{q^{w_{2}}}\left(y_{1}\right)\cdots d\mu_{q^{w_{2}}}\left(y_{l}\right)\nonumber \\
= & \sum_{i=0}^{n}\dbinom{n}{i}\left[w_{1}\right]_{q}^{n-i}\left[w_{2}\right]_{q}^{i-r}\beta_{i,q^{w_{2}}}^{\left(r\right)}\left(w_{1}x\right)T_{n,i}^{\left(r\right)}\left(w_{2}|q^{w_{1}}\right).\nonumber 
\end{align}

Therefore, by Corollary \ref{cor:2}, (\ref{eq:27}) and (\ref{eq:29}),
we obtain the following theorem.
\begin{thm}
\label{thm:4}For $n\ge0$ and $r,\, w_{1},\, w_{2}\in\mathbb{N}$,
we have 
\begin{align*}
 & \sum_{i=0}^{n}\dbinom{n}{i}\left[w_{1}\right]_{q}^{n-i}\left[w_{2}\right]_{q}^{i-r}\beta_{i,q^{w_{2}}}^{\left(r\right)}\left(w_{1}x\right)T_{n,i}^{\left(r\right)}\left(w_{2}|q^{w_{1}}\right)\\
= & \sum_{i=0}^{n}\dbinom{n}{i}\left[w_{2}\right]_{q}^{n-i}\left[w_{1}\right]_{q}^{i-r}\beta_{i,q^{w_{1}}}^{\left(r\right)}\left(w_{2}x\right)T_{n,i}^{\left(r\right)}\left(w_{1}|q^{w_{2}}\right),
\end{align*}
where 
\[
T_{n,i}^{\left(r\right)}\left(w|q\right)=\sum_{j_{1},\cdots,j_{r}=0}^{w-1}\left[j_{1}+\cdots+j_{r}\right]_{q}^{n-i}q^{\left(i+1\right)\sum_{l=1}^{r}j_{l}}.
\]

\end{thm}
$\,$

For $h\in\mathbb{Z}$ and $r\in\mathbb{N}$, we have 
\begin{align*}
 & \int_{\Zp}\cdots\int_{\Zp}\left[x+y_{1}+\cdots+y_{r}\right]_{q}^{n}q^{\sum_{l=1}^{r}\left(h-l\right)y_{l}}d\mu_{q}\left(y_{1}\right)\cdots d\mu_{q}\left(y_{r}\right)\\
= & \sum_{j=0}^{n}\dbinom{n}{j}\left(-1\right)^{j}\frac{q^{xj}}{\left(1-q\right)^{n}}\lim_{N\rightarrow\infty}\frac{1}{\left[p^{N}\right]_{q}^{r}}\sum_{y_{1},\cdots,y_{r}=0}^{p^{N}-1}q^{j\sum_{l=1}^{r}y_{l}}q^{\sum_{l=1}^{r}\left(h-l+1\right)y_{l}}\\
= & \sum_{j=0}^{n}\dbinom{n}{j}\left(-1\right)^{j}\frac{q^{xj}}{\left(1-q\right)^{n}}\frac{\left(j+h\right)\left(j+h-1\right)\cdots\left(j+h-r+1\right)}{\left[j+h\right]_{q}\left[j+h-1\right]_{q}\cdots\left[j+h-r+1\right]_{q}}\\
= & \frac{1}{\left(1-q\right)^{n}}\sum_{j=0}^{n}\dbinom{n}{j}\left(-1\right)^{j}q^{xj}\frac{\binom{j+h}{r}}{\binom{j+h}{r}_{q}}\frac{r!}{\left[r\right]_{q}!},
\end{align*}
where $\dbinom{x}{r}_{q}=\frac{\left[x\right]_{q}\left[x-1\right]_{q}\cdots\left[x-r+1\right]_{q}}{\left[r\right]_{q}!}=\frac{\left[x\right]_{q}\left[x-1\right]_{q}\cdots\left[x-r+1\right]_{q}}{\left[r\right]_{q}\left[r-1\right]_{q}\cdots\left[2\right]_{q}\left[1\right]_{q}}.$

From (\ref{eq:18}), we can also define $q$-extensions of higher-order
Bernoulli polynomials as follows : 

\begin{equation}
\beta_{n,q}^{\left(h,r\right)}\left(x\right)=\int_{\Zp}\cdots\int_{\Zp}\left[x+y_{1}+\cdots+y_{r}\right]_{q}^{n}q^{\sum_{l=1}^{r}\left(h-l\right)y_{l}}d\mu_{q}\left(y_{1}\right)\cdots d\mu_{q}\left(y_{r}\right),\label{eq:30}
\end{equation}
where $n\ge0$ and $h\in\mathbb{Z}$, $r\in\mathbb{N}$.

Let $w_{1},\, w_{2}\in\mathbb{N}$. Then we see that
\begin{align}
 & \frac{1}{\left[w_{1}\right]_{q}^{r}}\sum_{j_{1},\cdots,j_{r}=0}^{w_{1}-1}q^{w_{2}\sum_{l=1}^{r}\left(h-l+1\right)j_{l}}\label{eq:31}\\
 & \times\int_{\Zp}\cdots\int_{\Zp}q^{w_{1}\sum_{l=1}^{r}\left(h-l\right)y_{l}}e^{\left[w_{1}w_{2}x+\sum_{l=1}^{r}\left(w_{2}j_{l}+w_{1}y_{l}\right)\right]_{q}t}d\mu_{q^{w_{1}}}\left(y_{1}\right)\cdots d\mu_{q^{w_{1}}}\left(y_{r}\right)\nonumber \\
= & \frac{1}{\left[w_{2}\right]_{q}^{r}}\sum_{j_{1},\cdots,j_{r}=0}^{w_{2}-1}q^{w_{1}\sum_{l=1}^{r}\left(h-l+1\right)j_{l}}\nonumber \\
 & \times\int_{\Zp}\cdots\int_{\Zp}q^{w_{2}\sum_{l=1}^{r}\left(h-l\right)y_{l}}e^{\left[w_{1}w_{2}x+\sum_{l=1}^{r}\left(w_{1}j_{l}+w_{2}y_{l}\right)\right]_{q}t}d\mu_{q^{w_{2}}}\left(y_{1}\right)\cdots d\mu_{q^{w_{2}}}\left(y_{r}\right).\nonumber 
\end{align}

From (\ref{eq:31}), we have 
\begin{align}
 & \left[w_{1}\right]_{q}^{n-r}\sum_{j_{1},\cdots,j_{r}=0}^{w_{1}-1}q^{w_{2}\sum_{l=1}^{r}\left(h-l+1\right)j_{l}}\label{eq:32}\\
 & \times\int_{\Zp}\cdots\int_{\Zp}q^{w_{1}\sum_{l=1}^{r}\left(h-l\right)y_{l}}\left[w_{2}x+\frac{w_{2}}{w_{1}}\sum_{l=1}^{r}j_{l}+\sum_{l=1}^{r}y_{l}\right]_{q^{w_{1}}}^{n}d\mu_{q^{w_{1}}}\left(y_{1}\right)\cdots d\mu_{q^{w_{1}}}\left(y_{r}\right)\nonumber \\
= & \left[w_{2}\right]_{q}^{n-r}\sum_{j_{1},\cdots,j_{r}=0}^{w_{2}-1}q^{w_{1}\sum_{l=1}^{r}\left(h-l+1\right)j_{l}}\nonumber \\
 & \times\int_{\Zp}\cdots\int_{\Zp}q^{w_{2}\sum_{l=1}^{r}\left(h-l\right)y_{l}}\left[w_{1}x+\frac{w_{1}}{w_{2}}\sum_{l=1}^{r}j_{l}+\sum_{l=1}^{r}y_{l}\right]_{q^{w_{2}}}^{n}d\mu_{q^{w_{2}}}\left(y_{1}\right)\cdots d\mu_{q^{w_{2}}}\left(y_{r}\right),\nonumber 
\end{align}
where $n\ge0$ and $r\in\mathbb{N}$, $h\in\mathbb{Z}$.

Therefore, by (\ref{eq:30}) and (\ref{eq:32}), we obtain the following
theorem. 
\begin{thm}
\label{thm:5} For $n\ge0$, $h\in\mathbb{Z}$ and $w_{1},w_{2}\in\mathbb{N}$,
we have 
\begin{align*}
 & \left[w_{1}\right]_{q}^{n-r}\sum_{j_{1},\cdots,j_{r}=0}^{w_{1}-1}q^{w_{2}\sum_{l=1}^{r}\left(h-l+1\right)j_{l}}\beta_{n,q^{w_{1}}}^{\left(h,r\right)}\left(w_{2}x+\frac{w_{2}}{w_{1}}\left(j_{1}+\cdots+j_{r}\right)\right)\\
= & \left[w_{2}\right]_{q}^{n-r}\sum_{j_{1},\cdots,j_{r}=0}^{w_{2}-1}q^{w_{1}\sum_{l=1}^{r}\left(h-l+1\right)j_{l}}\beta_{n,q^{w_{2}}}^{\left(h,r\right)}\left(w_{1}x+\frac{w_{1}}{w_{2}}\left(j_{1}+\cdots+j_{r}\right)\right).
\end{align*}

\end{thm}
$\,$

From (\ref{eq:30}), we can derive the following equation : 
\begin{align}
 & \int_{\Zp}\cdots\int_{\Zp}q^{w_{1}\sum_{l=1}^{r}\left(h-l\right)y_{l}}\left[w_{2}x+\frac{w_{2}}{w_{1}}\sum_{l=1}^{r}j_{l}+\sum_{l=1}^{r}y_{l}\right]_{q^{w_{1}}}^{n}d\mu_{q^{w_{1}}}\left(y_{1}\right)\cdots d\mu_{q^{w_{1}}}\left(y_{r}\right)\label{eq:33}\\
= & \sum_{i=0}^{n}\dbinom{n}{i}\left(\frac{\left[w_{2}\right]_{q}}{\left[w_{1}\right]_{q}}\right)^{i}\left[j_{1}+\cdots+j_{r}\right]_{q^{w_{2}}}^{i}q^{w_{2}\left(n-i\right)\sum_{l=1}^{r}j_{l}}\nonumber \\
 & \times\int_{\Zp}\cdots\int_{\Zp}q^{w_{1}\sum_{l=1}^{r}\left(h-r\right)y_{l}}\left[w_{2}x+\sum_{l=1}^{r}y_{l}\right]_{q^{w_{1}}}^{n-i}d\mu_{q^{w_{1}}}\left(y_{1}\right)\cdots d\mu_{q^{w_{1}}}\left(y_{r}\right)\nonumber \\
= & \sum_{i=0}^{n}\dbinom{n}{i}\left(\frac{\left[w_{2}\right]_{q}}{\left[w_{1}\right]_{q}}\right)^{i}\left[j_{1}+\cdots+j_{r}\right]_{q^{w_{2}}}^{i}q^{w_{2}\left(n-i\right)\sum_{l=1}^{r}j_{l}}\beta_{n-i,q^{w_{1}}}^{\left(h,r\right)}\left(w_{2}x\right).\nonumber 
\end{align}

By (\ref{eq:33}), we get 
\begin{align}
 & \left[w_{1}\right]_{q}^{n-r}\sum_{j_{1},\cdots,j_{r}=0}^{w_{1}-1}q^{w_{2}\sum_{l=1}^{r}\left(h-l+1\right)j_{l}}\label{eq:34}\\
 & \times\int_{\Zp}\cdots\int_{\Zp}q^{w_{1}\sum_{l=1}^{r}\left(h-l\right)y_{l}}\left[w_{2}x+\frac{w_{2}}{w_{1}}\sum_{l=1}^{r}j_{l}+\sum_{l=1}^{r}y_{l}\right]_{q^{w_{1}}}^{n}d\mu_{q^{w_{1}}}\left(y_{1}\right)\cdots d\mu_{q^{w_{1}}}\left(y_{r}\right)\nonumber \\
= & \sum_{j_{1},\cdots,j_{r}=0}^{w_{1}-1}q^{w_{2}\sum_{l=1}^{r}\left(h-l+1\right)j_{l}}\sum_{i=0}^{n}\dbinom{n}{i}\left[w_{2}\right]_{q}^{i}\left[w_{1}\right]_{q}^{n-i-r}\left[j_{1}+\cdots+j_{r}\right]_{q^{w_{2}}}^{i}\nonumber \\
& \times q^{w_{2}\left(n-i\right)\sum_{l=1}^{r}j_{l}}\beta_{n-i,q^{w_{1}}}^{\left(h,r\right)}\left(w_{2}x\right)\nonumber \\
= & \sum_{i=0}^{n}\dbinom{n}{i}\left[w_{2}\right]_{q}^{i}\left[w_{1}\right]_{q}^{n-i-r}\beta_{n-i,q^{w_{1}}}^{\left(h,r\right)}\left(w_{2}x\right)\sum_{j_{1},\cdots,j_{r}=0}^{w_{1}-1} \left[j_{1}+\cdots+j_{r}\right]_{q^{w_{2}}}^{i}\nonumber \\
& \times q^{w_{2}\sum_{l=1}^{r}\left(n+h-l-i+1\right)j_{l}}\nonumber\\
= & \sum_{i=0}^{n}\dbinom{n}{i}\left[w_{2}\right]_{q}^{n-i}\left[w_{1}\right]_{q}^{i-r}\beta_{i,q^{w_{1}}}^{\left(h,r\right)}\left(w_{2}x\right)T_{n,i}^{\left(h,r\right)}\left(w_{1}|q^{w_{2}}\right),\nonumber 
\end{align}
where 
\begin{equation}
T_{n,i}^{\left(h,r\right)}\left(w|q\right)=\sum_{j_{1},\cdots,j_{r}=0}^{w-1}\left[j_{1}+\cdots+j_{r}\right]_{q}^{n-i}q^{\sum_{l=1}^{r}\left(i+h-l+1\right)j_{l}}.\label{eq:35}
\end{equation}

By the same method as (\ref{eq:34}), we see that 
\begin{align}
 & \left[w_{2}\right]_{q}^{n-r}\sum_{j_{1},\cdots,j_{r}=0}^{w_{2}-1}q^{w_{1}\sum_{l=1}^{r}\left(h-l+1\right)j_{l}}\label{eq:36}\\
 & \times\int_{\Zp}\cdots\int_{\Zp}q^{w_{2}\sum_{l=1}^{r}\left(h-l\right)y_{l}}\left[w_{1}x+\frac{w_{1}}{w_{2}}\sum_{l=1}^{r}j_{l}+\sum_{l=1}^{r}y_{l}\right]_{q^{w_{1}}}^{n}d\mu_{q^{w_{2}}}\left(y_{1}\right)\cdots d\mu_{q^{w_{2}}}\left(y_{r}\right)\nonumber \\
= & \sum_{i=0}^{n}\dbinom{n}{i}\left[w_{1}\right]_{q}^{n-i}\left[w_{2}\right]_{q}^{i-r}\beta_{i,q^{w_{2}}}^{\left(h,r\right)}\left(w_{1}x\right)T_{n,i}^{\left(h,r\right)}\left(w_{2}|q^{w_{1}}\right).\nonumber 
\end{align}

Therefore, by (\ref{eq:34}) and (\ref{eq:36}), we obtain the following
theorem.
\begin{thm}
\label{thm:6} For $n\ge0$, $h\in\mathbb{Z}$ and $r,\, w_{1},\, w_{2}\in\mathbb{N}$,
we have 
\begin{align*}
 & \sum_{i=0}^{n}\dbinom{n}{i}\left[w_{2}\right]_{q}^{n-i}\left[w_{1}\right]_{q}^{i-r}\beta_{i,q^{w_{1}}}^{\left(h,r\right)}\left(w_{2}x\right)T_{n,i}^{\left(h,r\right)}\left(w_{1}|q^{w_{2}}\right)\\
= & \sum_{i=0}^{n}\dbinom{n}{i}\left[w_{1}\right]_{q}^{n-i}\left[w_{2}\right]_{q}^{i-r}\beta_{i,q^{w_{2}}}^{\left(h,r\right)}\left(w_{1}x\right)T_{n,i}^{\left(h,r\right)}\left(w_{2}|q^{w_{1}}\right),
\end{align*}
where 
\[
T_{n,i}^{\left(h,r\right)}\left(w|q\right)=\sum_{j_{1},\cdots,j_{r}=0}^{w-1}\left[j_{1}+\cdots+j_{r}\right]_{q}^{n-i}q^{\sum_{l=1}^{r}\left(h+i-l+1\right)j_{l}}.
\]
\end{thm}
\begin{rem*}
A $p$-adic approach to identities of symmetry for Carlitz's $q$-Bernoulli
polynomials has been studied in \cite{key-10}.
\end{rem*}
$\,$

\textbf{Acknowledgement }

This work was supported by the National Research Foundation of Korea
(NRF) grant funded by the Korea government (MOE) (No.2012R1A1A2003786).

\bibliographystyle{amsplain}
\nocite{*}
\bibliography{0111}

$\,$

\noindent \textsc{Department of Mathematics, Sogang University, Seoul
121-742, Republic of Korea}

\noindent \emph{E-mail}\emph{address : }\texttt{dskim@sogang.ac.kr}

$\,$

\noindent \textsc{Department of Mathematics, Kwangwoon University, Seoul
139-701, Republic of Korea}

\noindent \emph{E-mail}\emph{address : }\texttt{tkkim@kw.ac.kr} 
\end{document}